\documentclass[12pt]{article}
\usepackage{graphicx}
\usepackage{color}
\catcode`\@=11 \@addtoreset{equation}{section}

\catcode`\@=12

\newtheorem{Theorem}{Theorem}[section]
\newtheorem{Proposition}{Proposition}[section]
\newtheorem{Lemma}{Lemma}[section]
\newtheorem{Corollary}{Corollary}[section]
\newtheorem{Remark}{Remark}[section]

\newcommand{\bTheorem}[1]{
\begin{Theorem} \label{T#1} }
\newcommand{\eT}{\end{Theorem}}

\newcommand{\bProposition}[1]{
\begin{Proposition} \label{P#1}}
\newcommand{\eP}{\end{Proposition}}

\newcommand{\bLemma}[1]{
\begin{Lemma} \label{L#1} }
\newcommand{\eL}{\end{Lemma}}

\newcommand{\bCorollary}[1]{
\begin{Corollary} \label{C#1} }
\newcommand{\eC}{\end{Corollary}}

\newcommand{\bRemark}[1]{
\begin{Remark} \label{R#1} }
\newcommand{\eR}{\end{Remark}}

\newcommand{\bFormula}[1]{
\begin{equation} \label{#1}}
\newcommand{\eF}{\end{equation}}

\newcommand{\Ov}[1]{\overline{#1}}

\newcommand{\DC}{C^\infty_c}

\newcommand{\vu}{\vc{u}}
\newcommand{\vc}[1]{{\bf #1}}

\newcommand{\Div}{{\rm div}_x}
\newcommand{\Grad}{\nabla_x}

\newcommand{\tn}[1]{\mbox {\F #1}}
\newcommand{\dx}{{\rm d} {x}}
\newcommand{\dt}{{\rm d} t }

\newcommand{\dxdt}{\dx \ \dt}
\newcommand{\intO}[1]{\int_{\Omega} #1 \ \dx}

\newcommand{\bProof}{{\bf Proof: }}

\font\F=msbm10 scaled 1000
\newcommand{\R}{\mbox{\F R}}

\newcommand{\Del}{\Delta_x}

\definecolor{Cgrey}{rgb}{0.85,0.85,0.85}
\definecolor{Cblue}{rgb}{0.50,0.85,0.85}
\definecolor{Cred}{rgb}{1,0,0}
\definecolor{fancy}{rgb}{0.10,0.85,0.10}

\newcommand\Cbox[2]{%
    \newbox\contentbox%
    \newbox\bkgdbox%
    \setbox\contentbox\hbox to \hsize{%
        \vtop{
            \kern\columnsep
            \hbox to \hsize{%
                \kern\columnsep%
                \advance\hsize by -2\columnsep%
                \setlength{\textwidth}{\hsize}%
                \vbox{
                    \parskip=\baselineskip
                    \parindent=0bp
                    #2
                }%
                \kern\columnsep%
            }%
            \kern\columnsep%
        }%
    }%
    \setbox\bkgdbox\vbox{
        \color{#1}
        \hrule width  \wd\contentbox %
               height \ht\contentbox %
               depth  \dp\contentbox
        \color{black}
    }%
    \wd\bkgdbox=0bp%
    \vbox{\hbox to \hsize{\box\bkgdbox\box\contentbox}}%
    \vskip\baselineskip%
}

\date{}

\begin{document}

%%%%%%%%%%%%%%%%%%%%%%%%%%%%%%%%

\title{On weak solutions to the 2D Savage-Hutter model of the motion of a gravity driven avalanche flow}

\author{Eduard Feireisl \thanks{The research of E.F. leading to these results has received funding from the European Research Council under the European Union's Seventh Framework
Programme (FP7/2007-2013)/ ERC Grant Agreement 320078} \and Piotr Gwiazda\thanks{P. G.
is a coordinator of the International Ph.D. Projects Programme of the Foundation for Polish Science operated
within the Innovative Economy Operational Programme 2007--2013 (Ph.D. Programme: Mathematical Methods
in Natural Sciences)} \and Agnieszka \' Swierczewska-Gwiazda\thanks{The research of A.\'S.-G. has received funding from the  National Science Centre, DEC-2012/05/E/ST1/02218} }

\maketitle

\bigskip

\centerline{Institute of Mathematics of the Academy of Sciences of the Czech Republic}

\centerline{\v Zitn\' a 25, 115 67 Praha 1, Czech Republic}

\bigskip

\centerline{Institute of Applied Mathematics and Mechanics, University of Warsaw}

\centerline{Banacha 2, 02-097 Warszawa, Poland}

%\centerline{Charles University in Prague, Faculty of Mathematics and Physics, Mathematical Institute}

%\centerline{Sokolovsk\' a 83, 186 75 Praha 8, Czech Republic}

%\thanks{Eduard Feireisl acknowledges the support of the GA\v CR (Czech Science Foundation) project 13-00522S
%in the general framework of RVO: 67985840}

%\thanks{Eduard Feireisl acknowledges the support of the project LL1202 in the
%programme ERC-CZ funded
%by the Ministry of Education, Youth and Sports of the Czech Republic.}

%\thanks{The research of E.F. leading to these results has received funding from the European Research Council under the European Union's Seventh %Framework
%Programme (FP7/2007-2013)/ ERC Grant Agreement 320078}

\bigskip

%\centerline{Institute of Mathematics of the Academy of Sciences of the Czech Republic}

%\centerline{\v Zitn\' a 25, 115 67 Praha 1, Czech Republic}

%\centerline{Charles University in Prague, Faculty of Mathematics and Physics, Mathematical Institute}

%\centerline{Sokolovsk\' a 83, 186 75 Praha 8,
%Czech Republic}

\begin{abstract}

We consider the Savage-Hutter system consisting of  two-dimensional depth-integrated shallow
water equations for the incompressible fluid with the
Coulomb-type friction term. Using the method of convex integration we show that the associated initial-value problem
possesses infinitely many weak solutions for any finite-energy initial data. On the other hand, the problem enjoys
the weak-strong uniqueness property provided the system of equations is supplemented with the energy inequality.

\end{abstract}


%\tableofcontents

\section{Introduction}
\label{i}

The theory for gravity driven avalanche flows is qualitatively similar to that of compressible fluid dynamics. We consider a relatively simple
Savage-Hutter model based on the assumption that the material is incompressible, with an isotropic pressure distribution through its depth, and a Coulomb
sliding friction, see Gray et al. \cite{GrTaNo} and \cite{HuWaPu} for a comprehensive overview. The time evolution of the flow height $h = h(t,x)$ and depth-averaged velocity $\vu = \vu(t,x)$ is
described through a system of balance laws - the \emph{Savage-Hutter system}:

\Cbox{Cgrey}{

\bFormula{i1}
\partial_t h + \Div ( h \vu ) = 0,
\eF
\bFormula{i2}
\partial_t (h \vu) + \Div (h \vu  \otimes \vu ) + \Grad ( a h^2 ) = h \left( - \gamma \frac{\vu}{ |\vu|} +  \vc{f} \right),
\eF

}

\noindent where $|\cdot|$ is an Euclidean metric,  $a \geq 0$, $\gamma \geq 0$, and $\vc{f}$ are given (smooth) functions of spatial coordinate $x \in \Omega \subset R^2$, see
Gray and Cui \cite{GraCui}, Zahibo et al. \cite{ZaPeTaNi}. For the sake of simplicity, we restrict ourselves to the periodic boundary conditions supposing accordingly that
$\Omega$ is the ``flat'' torus
\bFormula{i3}
\Omega = \left( [0,1]|_{ \{ 0, 1 \} } \right)^2.
\eF
The system (\ref{i1}), (\ref{i2}) is supplemented with the initial conditions
\bFormula{i4}
h(0, \cdot) = h_0, \ \vu(0, \cdot) = \vu_0.
\eF

The term $\frac{\vu}{ |\vu|}$ has to be  understood as a multi-valued mapping, which for non-zero velocities takes the value $\frac{\vu}{ |\vu|}$, whereas for $\vu=0$ takes the values in the whole closed unit ball. For the physical justification for such formulation with a simple argument on one-dimensional steady solutions we refer to \cite{Gw1}. 

In the absence of the driving force on the right-hand side of (\ref{i2}), the Savage-Hutter system coincides with barotropic Euler system describing the motion of a compressible inviscid fluid. As observed Gray and Cui \cite{GraCui}, the solutions of the Savage-Hutter system develop shock waves and other singularities characteristic for hyperbolic system of conservation laws. Accordingly, any mathematical theory based on the classical concept of (smooth) solutions fails
as soon as we are interested in global-in-time solutions to the system (\ref{i1}), (\ref{i2}), and/or in solutions emanating from singular initial data. 
Hence in analogue to the development of the theory for Euler equations, the problem of existence of weak solutions  has so far remained open. 
The issue of measure-valued solutions for the two-dimensional model considered here was studied in \cite{Gw2}.  Again, similar as for  Euler flows, the author follows the concept of generalization by DiPerna and Majda \cite{DiPeMa} to capture both oscillations as well as concentration effects. The existence of entropy weak solution to the corresponding problem in one-dimensional setting was shown in \cite{Gw1}.  Various modifications of the model to more complex topographies were considered in \cite{Boetal, BoWe, FeNietal}, see also \cite{JuMuGaNa, PeBoMa} for computational results. 

In this paper, we consider the \emph{weak solutions} to the problem (\ref{i1}--\ref{i4}) determined by means of a family of integral identities:

\bFormula{i5}
\int_0^T \intO{ \left( h \partial_t \varphi + h \vu \cdot \Grad \varphi \right) } \ \dt = - \intO{ h_0 \varphi (0, \cdot) } \ \dt
\eF
for any $\varphi \in C^1([0,T) \times \Omega)$;
\bFormula{i6}
\int_0^T \intO{ \left( h \vu \cdot \partial_t \Phi + h \vu \otimes \vu : \Grad \Phi + a h^2 \Div \Phi \right) } \ \dt
\eF
\[
= \int_0^T \intO{ h \left( \gamma \vc{B} - \vc{f} \right) \cdot \Phi } \ \dt - \intO{ h_0 \vu_0 \cdot \Phi (0, \cdot) }
\]
for any $\Phi \in C^1([0,T) \times \Omega; R^2)$,
where
\bFormula{i7}
\vc{B} = \vc{B}_{\vu} (t,x) = \left\{ \begin{array}{l} \frac{\vu(t,x)}{|\vu(t,x)|} \ \mbox{if}\ \vu(t,x) \ne 0, \\ \\ \in \overline{B_1(0)} \ \mbox{if}
\ \vu(t,x) = 0,      \end{array} \right.
\eF
where $B_1(0)$ is the unit ball in $\R^2$ with respect to the Euclidean metric. Hence speaking about weak solutions we in fact mean a triplet $[h,\vu,\vc{B}_{\vu}]$, whereas $\vc{B}_\vu$ is the selection from the multi-valued graph.

Using the method of convex integration, recently adapted to the incompressible Euler system by De Lellis and Sz\' ekelyhidi \cite{DelSze3}, and Chiodaroli \cite{Chiod},
we show that the Savage-Hutter system is always \emph{solvable} but not \emph{well posed} in the class of weak solutions. More specifically, we show
that the problem (\ref{i1}--\ref{i4}) admits \emph{infinitely many} weak solutions on a given time interval $(0,T)$ and for any sufficiently smooth initial data, see Section \ref{I}. 
The method of convex integration is used to construct large sets of weak solutions. Starting from the seminal results of Scheffer \cite{Sch} an Shnirelman \cite{Shn} on existence of compactly supported $2D$ weak solutions to the Euler equations and the previously mentioned  result of De Lellis and Sz\' ekelyhidi  \cite{DelSze3}. The theory was extensively developed for constructng  weak solutions with bounded energy by Wiedemann \cite{Wi} and in a nontrivial way   extended for compressible Euler in the so-called {\it variable coefficient} variant of the approach \cite{Chiod, ChKr} and recently to more complex systems such as equations of a compressible heat conducting gas \cite{ChiFeiKre}, Euler-Korteweg-Poisson system \cite{DoFeMa}.  

Apparently, the approach of convex integration allows to avoid the difficulties arising from the multi-valued formulation of the friction term. In case of measure-valued solutions which were generated by the approximate sequences of corresponding viscous problems one could not exclude the case of $\vu=0$, cf.~Gwiazda \cite{Gw2} and also the same difficulty in $1D$ case \cite{Gw1}. In the current framework we choose $h$ and the energy $E$ and  the consequent steps lead to finding the momentum. Although the friction term was nonlinear in terms of $\vu$, the situation  turns out to be significantly better as we express the friction term as
a product of a scalar function only of $h$ and $E$ and a function of $h$ and $hu$ which is linear.  The linearity then obviously provides weak continuity. And  what is the most essential, postulating $h>0$ and the energy being sufficiently large yields that the velocity $\vu$ is non-zero, and hence we are away from the set  where the friction term is multi-valued.

Next, in Section \ref{e}, we augment the weak formulation (\ref{i5}--\ref{i7}) by the energy inequality. Moreover, using the relative
energy method introduced by Dafermos \cite{Daf4} and developed in \cite{FeJiNo}, we show the weak-strong uniqueness principle: A weak solution satisfying the energy
inequality necessarily coincides with the strong solution emanating from the same initial data as long as the latter exists. The paper is concluded
by some remarks concerning the implications of the convex integration technique on the well-posedness problem considered in the class of finite energy
weak solutions, see Section \ref{c}.

We conclude the introduction with a remark that the current studies are a negative result for numerical simulation conducted for this system as the analysis shows that solutions with a sufficiently large energy exist but are non-unique. Nevertheless, this does not exclude there might exist solutions with a small energy, which are unique. 

\section{Infinitely many weak solutions}
\label{I}
We start with a brief description of the procedure of constructing the solutions. An essential tool is the Baire category method arising from the theory of differential inclusions. Note however that instead of
the abstract arguments one could provide the constructive, however essentially longer proof by adding oscillatory perturbations. 
Similarly to \cite{ChiFeiKre}, we will start with regular initial data 
\bFormula{I1}
h_0 \in C^2(\Omega), \ \vu_0 \in C^2(\Omega; R^2), h_0 > 0 \ \mbox{in} \ \Omega.
\eF
and reformulate the problem in different variables in order to obtain the system corresponding to an incompressible Euler system. The regularity of data appears to be important in the oscillatory lemma, which the same as in case of compressible Euler flow is formulated in the variable coefficients form with the coefficients generated by the data. 
For the new system we will construct a family of subsolutions. Indeed, we rewrite it as a linear system coupled with a nonlinear constraint.  Firstly, completing the space of 
subsolutions in appropriate topology we consider a family of functionals, which turn out to be lower-semicontinuous. In the next step we will conclude that as the limit object is a pointwise limit of continuous functions, namely it is \emph{Baire-1 function}, then the set of continuity points is infinite.
%Although we deal with the weak solution, the application of convex integration requires the initial data to be regular. Similarly to \cite{ChiFeiKre}, we suppose that
%\bFormula{I1}
%h_0 \in C^2(\Omega), \ \vu_0 \in C^2(\Omega; R^2), h_0 > 0 \ \mbox{in} \ \Omega.
%\eF
To show that the points where the functional vanishes are solutions to the system we will use an oscillatory lemma, namely Lemma~\ref{I1}. 

In the first step, using the standard Helmholtz decomposition, we may write
\bFormula{I2}
h_0 \vu_0 = \vc{v}_0 + \vc{V}_0 + \Grad \Psi_0, \ \mbox{where}\ \Div \vc{v}_0 = 0, \ \intO{ \Psi_0 } = 0,\ \intO{ \vc{v}_0 } = 0,\
\vc{V}_0 \in R^2.
\eF

The main result proved in the section reads:

\Cbox{Cgrey}{

\bTheorem{I1}

Let $T > 0$ and the initial data $h_0$, $\vu_0$ satisfying (\ref{I1}) be given. Suppose that $a, \gamma \in C^2(\Omega)$, $a, \gamma \geq 0$, $\vc{f} \in
C^1([0,T] \times \Omega; R^2)$.

Then the problem (\ref{i1}--\ref{i4}) admits infinitely many weak solutions in $(0,T) \times \Omega$. The weak solutions belong to the class
\[
h, \ \partial_t h, \Grad h \in C^1([0,T] \times \Omega),
\]
\[
\vu \in C_{\rm weak}([0,T]; L^2(\Omega;R^2)) \cap L^\infty((0,T) \times
\Omega; R^2),\ \Div \vu \in C([0,T] \times \Omega),
\]
\[
\vc{B}\in L^\infty((0,T) \times
\Omega; R^2).
\]

\eT

}
Above, by $C_{\rm weak}([0,T]; L^2(\Omega;R^2))$ we mean the space of functions continuous in time with respect to the weak topology of $L^2(\Omega;R^2)$. 
The remaining part of this section is devoted to the proof of Theorem \ref{TI1}.

\subsection{Convex integration ansatz}

Similarly to the decomposition (\ref{I2}), we look for solutions in the form
\[
h\vu = \vc{v} + \vc{V} + \Grad \Psi, \ \mbox{where}\ \Div \vc{v} = 0, \ \intO{ \Psi (t, \cdot) } = 0,\ \intO{ \vc{v}(t, \cdot) } = 0,\
\vc{V} = \vc{V}(t) \in R^2.
\]
Thus the continuity equation (\ref{i1}) reads
\bFormula{i3a}
\partial_t h + \Del \Psi = 0 \ \mbox{in}\ (0,T) \times \Omega, \ h(0, \cdot) = h_0, \ \Psi(0, \cdot) = \Psi_0.
\eF
Now, we can \emph{choose} $h = h(t,x) \in C^2([0,T] \times \Omega)$ in such a way that
\[
h(0, \cdot) = h_0,\ \partial_t h(0, \cdot) = - \Del \Psi_0, \ h(t, \cdot) > 0, \ \intO{ h(t, \cdot) } = \intO{ h_0 }\ \mbox{for all}\ t \in [0,T],
\]
and compute
\[
- \Del \Psi (t, \cdot) = \partial_t h(t, \cdot), \ \intO{ \Psi(t, \cdot) } = 0.
\]

Consequently, the original problem (\ref{i1}--\ref{i4}) reduces to finding the functions $\vc{v}$, $\vc{V}$ satisfying (weakly)
\bFormula{i7a}
\partial_t \vc{v} + \partial_t \vc{V} + \Div \left( \frac{ (\vc{v} + \vc{V} + \Grad \Psi ) \otimes (\vc{v} + \vc{V} + \Grad \Psi ) }{h}
+ \left( a h^2  + \partial_t \Psi \right) \tn{I} \right)
\eF
\[
= h \left( - \gamma \frac{ \vc{v} + \vc{V} + \Grad \Psi }{ | \vc{v} + \vc{V} + \Grad \Psi |} +  \vc{f} \right) ,
\]
\bFormula{i8a}
\Div \vc{v} = 0, \ \intO{ \vc{v}(t, \cdot) } = 0,
\eF
\bFormula{i9a}
\vc{v}(0, \cdot) = \vc{v}_0, \ \vc{V}(0) = \vc{V}_0.
\eF

\subsection{Kinetic energy}

We denote
\bFormula{i11}
E = \frac{1}{2} \frac{ | \vc{v} + \vc{V} + \Grad \Psi |^2 }{h}
\eF
the kinetic energy associated with the Savage-Hutter system. Analogously, we rewrite (\ref{i7a}) in the form
\bFormula{i7A}
\partial_t \vc{v} + \partial_t \vc{V} + \Div \left( \frac{ (\vc{v} + \vc{V} + \Grad \Psi ) \otimes (\vc{v} + \vc{V} + \Grad \Psi ) }{h}
- \frac{1}{2} \frac{ |\vc{v} + \vc{V} + \Grad \Psi |^2 }{h} \tn{I} \right)
\eF
\[
+ \Grad \left( E - \Lambda + a h^2  + \partial_t \Psi  \right)
= - \gamma \left( \frac{h}{2E} \right)^{1/2} \left( \vc{v} + \vc{V} + \Grad \Psi \right) +  h \vc{f} ,
\]
where $\Lambda = \Lambda(t)$ is a spatially homogeneous function to be determined below.
Finally, for
\bFormula{i11a}
E = \Lambda - a h^2 - \partial_t \Psi,
\eF
equation (\ref{i7A}) reduces to
\bFormula{i7B}
\partial_t \vc{v} + \partial_t \vc{V} + \Div \left( \frac{ (\vc{v} + \vc{V} + \Grad \Psi ) \otimes (\vc{v} + \vc{V} + \Grad \Psi ) }{h}
- \frac{1}{2} \frac{ |\vc{v} + \vc{V} + \Grad \Psi |^2 }{h} \tn{I} \right)
\eF
\[
= - \gamma \left( \frac{h}{2E} \right)^{1/2} \left( \vc{v} + \vc{V} + \Grad \Psi \right) +  h \vc{f} .
\]

\subsection{Determining $\vc{V}$}

The spatially homogeneous function $\vc{V}$ is determined as the unique solution of the ordinary differential equation
\bFormula{i7C}
\partial_t \vc{V} - \left[ \frac{1}{|\Omega| } \intO{ \gamma \left( \frac{h}{2E} \right)^{1/2} } \right] \vc{V} =
\frac{1}{|\Omega|} \intO{ \left[  \gamma \left( \frac{h}{2E} \right)^{1/2} \left( \vc{v} + \Grad \Psi \right) + h \vc{f} \right] },
\ \vc{V}(0) = \vc{V}_0.
\eF
Note that $\vc{V} = \vc{V}[\vc{v}]$ depends on $\vc{v}$ and also on the function $\Lambda$ in (\ref{i11a}).

With such a choice of $\vc{V}$, equation (\ref{i7B}) reads
\bFormula{i7D}
\partial_t \vc{v}  + \Div \left( \frac{ (\vc{v} + \vc{V}[\vc{v}] + \Grad \Psi ) \otimes (\vc{v} + \vc{V}[\vc{v}] + \Grad \Psi ) }{h}
- \frac{1}{2} \frac{ |\vc{v} + \vc{V}[\vc{v}] + \Grad \Psi |^2 }{h} \tn{I} \right)
\eF
\[
= - \gamma \left( \frac{h}{2E} \right)^{1/2} \left( \vc{v} + \vc{V}[\vc{v}] + \Grad \Psi \right)
+ \frac{1}{|\Omega|} \intO{ \gamma \left( \frac{h}{2E} \right)^{1/2} \left( \vc{v} + \vc{V}[\vc{v}] + \Grad \Psi \right) }
+  h \vc{f} - \frac{1}{|\Omega|} \intO{ h \vc{f} } .
\]

Finally, we find a tensor $\tn{M} = \tn{M}[\vc{v}]$ such that $\tn{M}(t,x) \in R^{2 \times 2}_{{\rm sym},0},$ where $R^{2 \times 2}_{{\rm sym},0}$ is the set of $2\times 2$ symmetric  traceless matrices,   for any $t,x$, and
\bFormula{I13}
\Div \tn{M}
\eF
\[
 =  - \gamma \left( \frac{h}{2E} \right)^{1/2} \left( \vc{v} + \vc{V}[\vc{v}] + \Grad \Psi \right)
+ \frac{1}{|\Omega|} \intO{ \gamma \left( \frac{h}{2E} \right)^{1/2} \left( \vc{v} + \vc{V}[\vc{v}] + \Grad \Psi \right) }
+  h \vc{f} - \frac{1}{|\Omega|} \intO{ h \vc{f} },
\]
\[
\intO{ \tn{M}(t, \cdot) } = 0 \ \mbox{for any}\ t \in [0,T].
\]
We can take
\[
\tn{M} = \Grad \vc{m} + \Grad^t \vc{m} - \Div \vc{m} \tn{I},
\]
where $\vc{m}$ is the (unique) solution of the elliptic equation
\bFormula{I14}
\Div \left( \Grad \vc{m} + \Grad^t \vc{m} - \Div \vc{m} \tn{I} \right)
\eF
\[
 =  - \gamma \left( \frac{h}{2E} \right)^{1/2} \left( \vc{v} + \vc{V}[\vc{v}] + \Grad \Psi \right)
+ \frac{1}{|\Omega|} \intO{ \gamma \left( \frac{h}{2E} \right)^{1/2} \left( \vc{v} + \vc{V}[\vc{v}] + \Grad \Psi \right) }
+  h \vc{f} - \frac{1}{|\Omega|} \intO{ h \vc{f} },
\]
\bFormula{I15}
\intO{ \vc{m} } = 0.
\eF
Indeed Desvilletes and Villani showed a variant of Korn's inequality \cite[Section IV.I, Proposition 11]{DesVil}
\[
\left\| \Grad \vc{m} + \Grad^t \vc{m} - \Div \vc{m} \tn{I} \right\|_{L^2(\Omega; R^{2 \times 2})} \geq \frac{1}{2} \| \Grad \vc{m} \|_{L^2(\Omega; R^{2 \times 2})},
\]
in particular, the problem (\ref{I14}), (\ref{I15}) admits a unique solution.

Summarizing, we write equation (\ref{i7D}) in a concise form
\bFormula{I16}
\partial_t \vc{v}  + \Div \left( \frac{ (\vc{v} + \vc{V}[\vc{v}] + \Grad \Psi ) \otimes (\vc{v} + \vc{V}[\vc{v}] + \Grad \Psi ) }{h}
- \frac{1}{2} \frac{ |\vc{v} + \vc{V}[\vc{v}] + \Grad \Psi |^2 }{h} \tn{I} - \tn{M}[\vc{v}] \right) = 0,
\eF
where $\vc{V}$, $\tn{M}$ are determined by means of (\ref{i7C}), (\ref{I13}), respectively.

\subsection{Application of the method of convex integration}

In order to recast our problem in terms of the method proposed by De Lellis and Sz\' ekelyhidi \cite{DelSze3}, we need to determine
\begin{itemize}
\item a (topological) space of subsolutions $X_0$;
\item a lower semi-continuous functional $I : \Ov{X}_0 \to R$ such that the points of continuity of $I$ coincide with the solution set of our problem.
\end{itemize}

\subsubsection{Subsolutions}

Let $\lambda_{\rm max}[\tn{A}]$ denote the maximal eigenvalue of a symmetric matrix $\tn{A}$.
Motivated by \cite{DelSze3}, we introduce the set
\[
X_0 = \left\{ \vc{w} \ \Big| \ {\vc{w}} \in L^\infty((0,T) \times \Omega, R^2) \cap C^1((0,T) \times \Omega,R^2) \cap
C_{\rm weak}([0,T]; L^2(\Omega;R^2)), \vphantom{\frac{1}{2}} \right.
\]
\[
\vc{w}(0, \cdot) = \vc{v}_0, \ \Div \vc{w} = 0 \ \mbox{in}\ (0,T) \times \Omega,
\]
\[
\partial_t \vc{w} + \Div \tn{F} = 0 \ \mbox{in}\ (0,T) \times \Omega
\ \mbox{for some}\ \tn{F} \in C^1((0,T) \times \Omega; R^{2 \times 2}_{0,{\rm sym}} ),
\]
\[
\lambda_{\rm max} \left[ \frac{ (\vc{w} + \vc{V}[\vc{w}] + \Grad \Psi) \otimes (\vc{w} + \vc{V}[\vc{w}] + \Grad \Psi) }{h} - \tn{F} - \tn{M}[\vc{w}] \right]
 < {E} - \delta
\]
\[
\left. \vphantom{\frac{1}{2}}
 \mbox{in}\ (0,T) \times \Omega \ \mbox{for some} \ \delta > 0 \right\},
\]
where $E$ is the kinetic energy introduced in (\ref{i11a}).

The first observation is that the set $X_0$ is non-empty provided
\bFormula{I17}
\Lambda (t) \geq \Lambda_0 > 0 \ \mbox{in} \ [0,T]
\eF
in (\ref{i11a}), and $\Lambda_0$ is large enough. Here ``large enough'' means in terms of the initial data, $\vc{f}$, and the time $T$.
Indeed, taking $\vc{w} = \vc{v_0}$, $\tn{F} = 0$, we have to find $\Lambda_0$ such that
\[
\lambda_{\rm max} \left[ \frac{ (\vc{v}_0 + \vc{V}[\vc{v}_0] + \Grad \Psi) \otimes (\vc{v}_0 + \vc{V}[\vc{v}_0] + \Grad \Psi) }{h}  - \tn{M}[\vc{v}_0] \right]
 < {E} - \delta = \Lambda - a h^2 - \partial_t \Psi - \delta.
\]
Since $\vc{V}$ is given by (\ref{i7C}), it is easy to check that $\vc{V}[\vc{v}_0]$, together with $\partial_t \vc{V}[\vc{v}_0]$, remain bounded in terms of the
data and uniformly for all $\Lambda \geq \Lambda_0$. Furthermore, applying the standard elliptic estimates to (\ref{I13}), we get
\[
\| \tn{M}(t, \cdot) \|_{W^{1,q}(\Omega; R^{2 \times 2})} \leq c(q, {\rm data}) \ \mbox{for any} \ 1 < q < \infty,\ t \in [0,T].
\]
Consequently, we may fix $\Lambda_0$, $\Lambda$ satisfying (\ref{I17}), and, finally, the kinetic energy $E$ in (\ref{i7A}) in such a way that the set of
subsolutions $X_0$ is non-empty.

\subsubsection{Uniform bounds}

As shown by De Lellis and Sz\' ekelyhidi \cite{DelSze3}, we have
\bFormula{I18}
\frac{1}{2} \frac{ |\vc{w} + \vc{V}[\vc{w}]  + \Grad \Psi |^2 }{h}
\eF
\[
\leq \lambda_{\rm max} \left[ \frac{ (\vc{w} + \vc{V}[\vc{w}] + \Grad \Psi) \otimes (\vc{w} + \vc{V}[\vc{w}] + \Grad \Psi) }{h} - \tn{F} - \tn{M}[\vc{w}] \right],
\]
where the equality holds only if
\bFormula{I19}
\tn{F} + \tn{M}[\vc{w}] = \frac{ (\vc{w} + \vc{V}[\vc{w}] + \Grad \Psi) \otimes (\vc{w} + \vc{V}[\vc{w}] + \Grad \Psi) }{h} -
\frac{ |\vc{w} + \vc{V}[\vc{w}] + \Grad \Psi |^2 }{h} \tn{I}.
\eF
Since $E$ has been fixed, we may deduce from (\ref{I18}) that
\bFormula{I19bis}
| \vc{w} + \vc{V}[\vc{w}] | \leq c(E,T, {\rm data}) \ \mbox{in}\ (0,T) \times \Omega,
\eF
and, going back to (\ref{i7C}), we may infer that
\bFormula{I20}
|\vc{w}| + |\vc{V}[\vc{w}]| + |\partial_t \vc{V}[\vc{w}] | \leq c(E, T, {\rm data}) \ \mbox{in}\ (0,T) \times \Omega.
\eF
Furthermore, (\ref{I14}) yields
\bFormula{I21}
\| \tn{M}[\vc{w}](t, \cdot) \|_{W^{1,q}(\Omega; R^{2 \times 2})} \leq c(E, T, {\rm data}) \ \mbox{in}\ (0,T) \ \mbox{for any}\ 1 \leq q < \infty.
\eF

The set $X_0$ is endowed with the topology of the space $C_{\rm weak}([0,T]; L^2(\Omega; R^2))$. In view of (\ref{I19bis}), such a topology is metrizable on
$X_0$ and we denote $\Ov{X}_0$ the completion of $X_0$ - a topological metric space. In accordance with (\ref{I20}), (\ref{I21}), and the compact embedding
$W^{1,q} \hookrightarrow\hookrightarrow C(\Omega)$, $q > 2$, we obtain
\bFormula{I22}
\left.
\begin{array}{r}
\vc{V}[\vc{w}_n] \to \vc{V}[\vc{w}] \ \mbox{in} \ C[0,T] \\ \\ \tn{M}[\vc{w}_n] \to \tn{M}[\vc{w}]
\ \mbox{in}\ C([0,T] \times \Omega) \end{array} \right\} \ \mbox{whenever}\ \vc{w}_n \to \vc{w} \ \mbox{in}\
C_{\rm weak}([0,T]; L^2(\Omega; R^2)).
\eF

\subsubsection{Functional $I$ and infinitely many solutions}

Following De Lellis and Sz\' ekelyhidi \cite{DelSze3}, we introduce the functional
\[
I[\vc{v}] = \int_0^T \intO{ \left[ \frac{1}{2} \frac{ |\vc{v} + \vc{V}[\vc{v}] + \Grad \Psi |^2 }{h} - E \right] }\ \dt: \Ov{X}_0 \to R.
\]
In order to proceed, we need the following variant of the oscillatory lemma (cf. De Lellis
and Sz\' ekelyhidi \cite[Proposition 3]{DelSze3}, Chiodaroli \cite[Section 6, formula (6.9)]{Chiod}) proved in \cite[Lemma 3.1]{DoFeMa} :

\bLemma{I1}

Let $U \subset R \times R^N$, $N=2,3$ be a bounded open set. Suppose that
\[
\vc{g} \in C(U; R^N), \ \tn{W} \in C(U; R^{N \times N}_{{\rm sym},0}), \ e,r \in C(U),\ r > 0, \ e \leq \Ov{e} \ \mbox{in}\ U
\]
are given such that
\[
\frac{N}{2} \lambda_{\rm max} \left[ \frac{ \vc{g} \otimes \vc{g} }{r} - \tn{W} \right] < e  \ \mbox{in}\
U.
\]
Then there exist sequences
\[
\vc{w}_n \in \DC (U;R^N), \ \tn{G}_n \in \DC(U; R^{N \times N}_{\rm sym,0}), \ n = 0,1,\dots
\]
such that
\[
\partial_t \vc{w}_n + \Div \tn{G}_n = 0 , \ \Div \vc{w}_n = 0 \ \mbox{in}\ R^N,
\]
\[
\frac{N}{2} \lambda_{\rm max} \left[ \frac{ (\vc{g} + \vc{w}_n) \otimes (\vc{g} + \vc{w}_n)}{r} - (\tn{W} + \tn{G}_n) \right] < e  \ \mbox{in}\ U,
\]
and
\bFormula{I23}
\vc{w}_n \to 0 \ \mbox{weakly in}\ L^2(U; R^N),\
\liminf_{n \to \infty} \int_{U} \frac{ | \vc{w}_n |^2 }{r} \ \dxdt \geq c(\Ov{e}) \int_{U} \left( e - \frac{1}{2} \frac{|\vc{g} |^2}{r} \right)^2 \ \dxdt.
\eF

\eL

\bRemark{I1}

It is important to note that the constant $c(\Ov{e})$ in (\ref{I23}) is independent of the specific form of the quantities $\vc{g}$, $\tn{W}$, $e$, and $r$.

\eR

In view of (\ref{I18}), we have
\[
I[\vc{w}] < 0 \ \mbox{for any}\ \vc{w} \in X_0,
\]
and, as a consequence of (\ref{I22}), $I: \Ov{X}_0 \to (- \infty, 0]$ is a lower semi-continuous functional with respect to the topology of
the space $C_{\rm weak}([0,T]; L^2(\Omega; R^2))$. Consequently, by virtue of Baire's category argument, the set of points
of continuity of $I$ in $\Ov{X}_0$ has infinite cardinality. Our ultimate goal will be to show that
\bFormula{I24}
I[\vc{v}] = 0 \ \mbox{whenever}\ \vc{v} \in \Ov{X}_0 \ \mbox{is a point of continuity of} \ I \ \mbox{in} \ \Ov{X}_0.
\eF
In view of (\ref{i11a}), (\ref{i7C}), (\ref{I14}), and (\ref{I18}), (\ref{I19}), it is easy to check that $\vc{v}$ represents a weak solution
of the problem (\ref{i7a}--\ref{i9a}), which completes the proof of Theorem \ref{TI1}.

To see (\ref{I24}),
arguing by contradiction, we assume that $\vc{v} \in \Ov{X}_0$ is a point of continuity of $I$ such that
\[
I[\vc{v}] < 0.
\]
Since $I$ is continuous at $\vc{v}$, there exists a sequence $\{ \vc{v}_m \}_{m=1}^\infty \subset X_0$ (and the associated fluxes $\tn{F}_m$) such that
\[
\vc{v}_m \to \vc{v} \ \mbox{in}\ C_{\rm weak}([0,T]; L^2(\Omega;R^2)), \ I[\vc{v}_m] \to I[\vc{v}] \ \mbox{as}\ m \to \infty.
\]
As $[\vc{v}_m, \tn{F}_m]$ are subsolutions, we get
\[
\lambda_{\rm max} \left[ \frac{ (\vc{v}_m + \vc{V}[\vc{v}_m] + \Grad \Psi) \otimes (\vc{v}_m + \vc{V}[\vc{v}_m] + \Grad \Psi) }{h} - \tn{F}_m
- \tn{M}[\vc{v}_m] \right]
\]
\[
< {E} - \delta_m \ \mbox{for some}\ \delta_m \searrow  0.
\]
Now,
fixing $m$ for a while, we apply Lemma \ref{LI1} with
\[
N=2, \ U = (0,T) \times \Omega,\ r = h, \ \vc{g} = \vc{v}_m + \vc{V}[\vc{v}_m] + \Grad \Psi, \ \tn{W} = \tn{F}_m + \tn{M}[\vc{v}_m],\
\mbox{and}\ e = E - \delta_m/2.
\]
Denoting $\left\{ [\vc{w}_{m,n}, \tn{G}_{m,n}] \right\}_{n=1}^\infty$ the quantities resulting from the conclusion of Lemma \ref{LI1}, we consider
\[
\vc{v}_{m,n} = \vc{v}_m + \vc{w}_{m,n}, \ \tn{F}_{m,n} = \tn{F}_m + \tn{G}_{m,n}.
\]
Obviously,
\[
\partial_t \vc{v}_{m,n} + \Div \tn{F}_{m,n} = 0,\ \Div \vc{v}_{m,n} = 0,\ \vc{v}_{m,n} (0, \cdot) = \vc{v}_0,
\]
and, in accordance with Lemma \ref{LI1},
\[
\lambda_{\rm max} \left[ \frac{ (\vc{v}_{m,n} + \vc{V}[\vc{v}_m] + \Grad \Psi) \otimes (\vc{v}_{m,n} + \vc{V}[\vc{v}_m] + \Grad \Psi) }{h} - \tn{F}_{m,n}
- \tn{M}[\vc{v}_m] \right]
< {E} - \delta_m/2.
\]
Consequently, in view of the continuity properties of the operators $\vc{v} \mapsto \vc{V}[\vc{v}]$, $\vc{v} \mapsto \tn{M} [\vc{v}]$, we may conclude that
for each $m$ there exists $n=n(m)$ such that
\[
[ \vc{v}_{m,n(m)}, \tn{F}_{m,n(m)} ] \in X_0, \ m=1,2,\dots
\]
Moreover, in view of (\ref{I23}), we may suppose
\[
\vc{v}_{m,n(m)} \to \vc{v} \ \mbox{in}\ C_{\rm weak}([0,T]; L^2(\Omega;R^2)) \ \mbox{as}\ m \to \infty,
\]
in particular,
\bFormula{I25}
I[\vc{v}_{m,n} ] \to I[\vc{v}] \ \mbox{as}\ m \to \infty.
\eF
Since for each $m$
\[
\lim_{n \to \infty} \int_0^T \intO{ \frac{( \vc{v}_m + \vc{V}[\vc{v}_m + \vc{w}_{m,n} ] + \Grad \Psi)\  \vc{w}_{m,n}}{h}  } \ \dt =0,
\]
thus 
using once more Lemma \ref{LI1} combined with Jensen's inequality, we observe that the sequence $\vc{v}_{m,n(m)}$ can be taken in such a way that
\[
\liminf_{m \to \infty} I [\vc{v}_{m,n(m)}] = \liminf_{m \to \infty} \int_0^T \intO{ \left( \frac{1}{2} \frac{ |\vc{v}_m + \vc{w}_{m,n(m)} + 
\vc{V}[\vc{v}_m + \vc{w}_{m,n(m)}] + \Grad \Psi|^2 }{h} - {E} \right) } \ \dt
\]
\[
= \lim_{m \to \infty} \int_0^T \intO{ \left( \frac{1}{2} \frac{ |\vc{v}_m + \vc{V}[\vc{v}_m + \vc{w}_{m,n(m)} ] + \Grad \Psi|^2 }{h} - {E} \right) } \ \dt + \liminf_{m \to \infty} \int_0^T \intO{ \frac{1}{2} \frac{ |\vc{w}_{m,n(m)}|^2 }{h} } \ \dt
\]
\[
\geq I [\vc{v}] + C_1 \liminf_{m \to \infty} \int_0^T \intO{ \left( E -\delta_m - \frac{1}{2} \frac{|\vc{v}_m + \vc{V}[\vc{v}_m] + \Grad \Psi|^2 }{h} \right)^2 } \ \dt 
\]
\[
\geq
I[\vc{v}] + C_2(T, |\Omega|) \liminf_{m \to \infty} \left( \int_0^T \intO{ \left( E -\delta_m - \frac{1}{2} \frac{|\vc{v}_m + \vc{V}[\vc{v}_m] + \Grad \Psi|^2 }{h} \right) } \ \dt   \right)^2
\]
\[
=
I[\vc{v}] + C_2(T, |\Omega|) \left( I[\vc{v}] \right)^2,\ C_2(T, |\Omega|) > 0,  
\]
which is compatible with (\ref{I25}) only if $I[\vc{v}] = 0$.

Thus we have shown (\ref{I24}), and, consequently, Theorem \ref{TI1}. 

\section{Dissipative solutions}
\label{e}

The solutions ``constructed'' in the proof of Theorem \ref{TI1} satisfy (\ref{i11a}), more specifically, 
\[
\frac{1}{2} h |\vu|^2 = E = \Lambda - a h^2 - \partial_t \Psi \ \mbox{for a.a.}\ (t, x) \in (0,T) \times \Omega.
\]
In particular, as $\Lambda$ has been chosen large, the total energy $E_{\rm tot}$ of the flow,
\[
E_{\rm tot}(t) = \intO{ \left[ \frac{1}{2} h |\vu |^2 + a h^2 \right](t, \cdot) } 
\]
may (and does in ``most'' cases) experience a jump at the initial time, 
\bFormula{e1}
\liminf_{t \to 0+} E_{\rm tot}(t) > \intO{ \left[ \frac{1}{2} h_0 |\vu_0 |^2 + a h_0^2 \right]}.
\eF
Apparently, solutions satisfying (\ref{e1}) are ``non-physical'' violating the First law of thermodynamics, at least if the forces $\vc{f}$ are regular. This observation leads to a natural admissibility criterion 
based on the energy balance appended to the definition of weak solutions to eliminate the oscillatory solutions constructed in Theorem 
\ref{TI1}. 

\subsection{Energy inequality and dissipative solutions} 
\label{ei}

For the sake of simplicity, suppose that $a > 0$ is a positive constant independent of $x$. Taking, formally, the scalar product of equation (\ref{i2}) with $\vu$ and integrating the resulting expression over $\Omega\times (0,\tau)$, we obtain the energy inequality 
\bFormula{e2}
E_{\rm tot}(\tau) \equiv \intO{ \left[ \frac{1}{2} h |\vu|^2 + a h^2 \right](\tau, \cdot) } + 
\int_0^\tau \intO{ h\gamma \vc{B}_{\vu} \cdot \vu } \ \dt \leq \intO{ \left[ \frac{1}{2} h_0 |\vu_0|^2 + a h^2_0 \right] } 
\eF
\[
+ \int_0^\tau \intO{ h \vc{f} \cdot \vu } \ \dt,
\] 
where the function $\vc{B}$ was introduced in (\ref{i7}).

We say that $[h, \vu,\vc{B}_\vu]$ is a \emph{dissipative} weak solution to the Savage-Hunter system if, in addition to (\ref{i5}--\ref{i7}), the energy 
inequality (\ref{e2}) holds for a.a. $\tau \in (0,T)$.  

\subsection{Relative energy and weak-strong uniqueness}

Our goal is to show that a dissipative and a strong solution emanating from the same initial data coincide as long as the latter exists. To this end, we revoke the method proposed by Dafermos \cite{Daf4} and later elaborated in \cite{FeJiNo}, based on the concept of \emph{relative energy}. We introduce
the relative energy functional
\bFormula{e3}
\mathcal{E} \left( h , \vu \ \Big|\ H, \vc{U} \right) = \intO{ \left[ \frac{1}{2} h |\vu - \vc{U}|^2 + P(h) - P'(H)(h - H) - P(H) \right] }, 
\eF
where 
\[
P(h) = a h^2.
\]

Now, exactly as in \cite[Section 3]{FeJiNo}, we may derive the relative energy inequality
\bFormula{e4}
{\mathcal E} \left( h , \vu \ \Big|\ H, \vc{U} \right)  (\tau) + \int_0^\tau \intO{   h\gamma \vc{B}_{\vu} \cdot (\vu - \vc{U})  } \ \dt \leq  {\mathcal E} \left( h_0 , \vu_0 \ \Big|\ H(0, \cdot), \vc{U}(0, \cdot) \right)
\eF
\[
+
\int_0^\tau \intO{  h \Big( \partial_t \vc{U} + \vu \Grad \vc{U} \Big) \cdot
(\vc{U} - \vu )} \ \dt 
+ \int_0^\tau \intO{ h \vc{f} \cdot(\vc u-\vc{U})} \ \dt
\]
\[
+ \int_0^\tau \intO{ \left( (H - h) \partial_t P'(H) + \Grad P'(H) \cdot
\left( H \vc{U} - h \vu \right) \right) } \ \dt
- \int_0^\tau \intO{
\Div \vc{U} \Big( a h^2 - a H^2 \Big) },
\] 
where $[h, \vu, \vc{B}_\vu]$ is a dissipative weak solution of the Savage-Hutter system and $[H, \vc{U},\vc{B}_\vc{U}]$ are sufficiently smooth functions in 
$[0,T] \times \Omega$, $H > 0$. 

We are ready to prove the following result:

\Cbox{Cgrey}{

\bTheorem{e1}

Let $[h,\vu,  \vc{B}_\vu]$ be a dissipative weak solution of the Savage-Hutter system in $(0,T) \times \Omega$  in the sense specified through (\ref{i5}--\ref{i7}),
(\ref{e2}). Let $[H, \vc{U}, \vc{B}_\vc{U}]$, $H > 0$ be a globally Lipschitz (strong) solution of the same problem, with 
\[
h_0 = H(0, \cdot),\ 
\vu_0 = \vc{U}(0, \cdot).
\]

Then 
\[
h = H,\ \vu = \vc{U} \ \mbox{a.e. in}\ (0,T) \times \Omega.
\]

\eT

}

As $h,\vu$ are uniquely determined, then from the balance of momentum one can recover  $
\vc{B}_\vc{u}$ such that $\vc{B}_\vc{u}=\vc{B}_\vc{U}$ for almost all $(t,x)\in (0,T) \times \Omega$.

\bProof

A simple density arguments shows that $[H, \vc{U}]$ may be used as test functions in the relative energy inequality (\ref{e4}). As the 
initial values coincide, the latter reads 
\bFormula{e5}
{\mathcal E} \left( h , \vu \ \Big|\ H, \vc{U} \right)  (\tau) + \int_0^\tau \intO{   h\gamma \vc{B}_{\vu} \cdot (\vu - \vc{U})  } \ \dt 
\eF
\[
\leq
\int_0^\tau \intO{  h \Big( \partial_t \vc{U} + \vu \Grad \vc{U} \Big) \cdot
(\vc{U} - \vu )} \ \dt
+ \int_0^\tau \intO{ h \vc{f} \cdot(\vc u-\vc{U})} \ \dt
\]
\[
+ 2a \int_0^\tau \intO{ \left( (H - h) \partial_t H + \Grad H \cdot
\left( H \vc{U} - h \vu \right) \right) } \ \dt
- \int_0^\tau \intO{
\Div \vc{U} \Big( a h^2 - a H^2 \Big) },
\] 
where, furthermore, 
\bFormula{e6}
h \Big( \partial_t \vc{U} + \vu \cdot \Grad \vc{U} \Big)\cdot (\vc{U} - \vu ) 
\eF
\[
= h \Big( \partial_t \vc{U} + \vc{U} \cdot \Grad \vc{U} \Big) \cdot (\vc{U} - \vu )
+ h (\vu - \vc{U}) \cdot \Grad \vc{U} \cdot (\vc{U} - \vu).
\]

As $\vc{U}$ is globally Lipschitz, the second term in (\ref{e6}) may be ``absorbed'' by the left-hand side of 
(\ref{e5}) via Gronwall's argument. Furthermore, 
\[
h \Big( \partial_t \vc{U} + \vc{U} \cdot \Grad \vc{U} \Big) \cdot (\vc{U} - \vu ) = 
\frac{h}{H} \Big( \partial_t (H \vc{U}) + \Div (H \vc{U} \otimes \vc{U}) \Big)\cdot (\vc{U} - \vu ) 
\]
\[
= - h\gamma \vc{B}_{\vc{U}} \cdot (\vc{U} - \vu ) - h \vc{f} \cdot (\vu - \vc{U}) - \frac{h}{H} \Grad (a H^2 ) \cdot (\vc{U} - \vu)  
\]
Consequently, (\ref{e5}) reduces to 
\bFormula{e7}
{\mathcal E} \left( h , \vu \ \Big|\ H, \vc{U} \right)  (\tau) + \int_0^\tau \intO{   h\gamma ( \vc{B}_{\vu} - \vc{B}_{\vc{U} } ) \cdot (\vu - \vc{U})  } \ \dt
\eF
\[
\leq c \int_0^\tau{\mathcal E} \left( h , \vu \ \Big|\ H, \vc{U} \right)  (t)
  \dt 
\]
\[
+ 2a \int_0^\tau \intO{ \left( (H - h) \partial_t H + \Grad H \cdot
\left( H \vc{U} - h \vc{U} \right) \right) } \ \dt
- \int_0^\tau \intO{
\Div \vc{U} \Big( a h^2 - a H^2 \Big) }.
\] 

Finally, write 
\[
\Div \vc{U} \Big( a h^2 - a H^2 \Big) = \Div \vc{U} \Big( a h^2 - 2a H (h - H) - a H^2 \Big) + 2a \Div \vc{U} H (h - H),
\]
where, similarly to the above, the first term can be handled by Gronwall's argument. Seeing the $H$, $\vc{U}$ satisfy the equation of continuity 
(\ref{i1}) we conclude that 
\[
{\mathcal E} \left( h , \vu \ \Big|\ H, \vc{U} \right)  (\tau) + \int_0^\tau \intO{   h\gamma ( \vc{B}_{\vu} - \vc{B} _{\vc{U} } ) \cdot (\vu - \vc{U})  } \ \dt
\leq c \int_0^\tau{\mathcal E} \left( h , \vu \ \Big|\ H, \vc{U} \right)  (t)
  \dt;
\]
whence, by Gronwall's lemma, ${\mathcal E} \left( h , \vu \ \Big|\ H, \vc{U} \right)(\tau) = 0$ for a.a. $\tau \in (0,T)$. 

\rightline{Q.E.D.}
 
\section{Concluding remarks}
\label{c}

In view of Theorem \ref{Te1}, we are tempted to say that imposing the energy inequality (\ref{e2}) eliminates the ``non-physical'' oscillatory 
solutions, the existence of which is claimed in Theorem \ref{TI1}. However, the method of convex integration may be used to obtain the following result that can be shown in the same way as \cite[Theorem 4.2]{ChiFeiKre}. 

{\bf Claim.} {\it Under the hypotheses of Theorem \ref{TI1}, let $T > 0$, $h_0$ be given. 
Then there exists 
\[
\vu_0 \in L^\infty(\Omega; R^2) 
\]
such that the Savage-Hutter system admits infinitely many dissipative weak solutions in $(0,T) \times \Omega$ starting from the initial data 
$[h_0, \vu_0]$.
}

\medskip
We leave the proof of {\bf Claim} to the interested reader. 

Finally we would like to mention that the same results hold in case of different choice of friction 
term. According to \cite{HuWaPu} the Savage-Hutter model is well valid for sand avalanches, however for the case of snow avalanches there is often considered a second velocity-dependent contribution, e.g. $h|\vu|\vu$. In a consequence, in terms of new variables the friction term reads for some coefficients $\gamma_1,\gamma_2\ge0$ as follows
$- \left( \gamma_1 \left( \frac{h}{2E} \right)^{1/2} +  \gamma_2 \left( \frac{2E}{h} \right)^{1/2} \right)\left( \vc{v} + \vc{V} + \Grad \Psi \right) $.

\def\cprime{$'$} \def\ocirc#1{\ifmmode\setbox0=\hbox{$#1$}\dimen0=\ht0
  \advance\dimen0 by1pt\rlap{\hbox to\wd0{\hss\raise\dimen0
  \hbox{\hskip.2em$\scriptscriptstyle\circ$}\hss}}#1\else {\accent"17 #1}\fi}

%\bibliography{citace}

\begin{thebibliography}{10}

\bibitem{Boetal}
 F.~Bouchut, A.~Mangeney-Castelnau, B.~Perthame, J.-P.~Vilotte.
 \newblock A new model of Saint Venant and Savage-Hutter type for gravity driven shallow water flows. 
 \newblock {\em C. R. Math. Acad. Sci. Paris} {\bf  336}(6):531--536, 2003.


\bibitem{BoWe}
F.~Bouchut, M.~Westdickenberg.
\newblock Gravity driven shallow water models for arbitrary topography. 
\newblock {\em Commun. Math. Sci.} {\bf 2}(3):359--389, 2004.

\bibitem{Chiod}
E.~Chiodaroli.
\newblock A counterexample to well-posedness of entropy solutions to the
  compressible {E}uler system.
\newblock {\em J. Hyperbolic Differ. Equ.}, {\bf 11}(3):493--519, 2014.

\bibitem{ChiFeiKre}
E.~Chiodaroli, E.~Feireisl, and O.~Kreml.
\newblock On the weak solutions to the equations of a compressible heat
  conducting gas.
\newblock {\em  Ann. Inst. H. Poincar\'e Anal. Non Lin\'eaire} {\bf 32}(1): 225--243, 2015.

\bibitem{ChKr} 
E.~Chiodaroli, O.~Kreml.
\newblock On the energy dissipation rate of solutions to the compressible isentropic Euler system. \newblock {\em Arch. Ration. Mech. Anal. } {\bf 214}(3): 1019--1049, 2014.

\bibitem{Daf4}
C.M. Dafermos.
\newblock The second law of thermodynamics and stability.
\newblock {\em Arch. Rational Mech. Anal.}, {\bf 70}:167--179, 1979.

\bibitem{DelSze3}
C.~De~Lellis and L.~Sz{\'e}kelyhidi, Jr.
\newblock On admissibility criteria for weak solutions of the {E}uler
  equations.
\newblock {\em Arch. Ration. Mech. Anal.}, {\bf 195}(1):225--260, 2010.

\bibitem{DesVil}
L.~Desvillettes and C.~Villani.
\newblock On the trend to global equilibrium for spatially inhomogeneous
  kinetic systems: the {B}oltzmann equation.
\newblock {\em Invent. Math.}, {\bf 159}(2):245--316, 2005.

\bibitem{DiPeMa}
 R.J.~DiPerna, A.J.~Majda, 
 \newblock Oscillations and concentrations in weak solutions of the incompressible fluid equations. \newblock {\em Comm. Math. Phys.} {\bf 108}(4): 667--689, 1987.

\bibitem{DoFeMa}
D.~Donatelli, E.~Feireisl, and P.~Marcati.
\newblock Well/ill posedness for the {E}uler-{K}orteweg-{P}oisson system and
  related problems.
\newblock {\em Commun. Partial Differential Equations}, 2014.
\newblock To appear.

\bibitem{FeJiNo}
E.~Feireisl, Bum~Ja Jin, and A.~Novotn{\' y}.
\newblock Relative entropies, suitable weak solutions, and weak-strong
  uniqueness for the compressible {N}avier-{S}tokes system.
\newblock {\em J. Math. Fluid Mech.}, {\bf 14}:712--730, 2012.

\bibitem{FeNietal}
 \newblock E.D.~Fern\'andez-Nieto, F.~Bouchut, D.~Bresch, M.J.~Castro D'az, A.~Mangeney.
 \newblock A new Savage-Hutter type model for submarine avalanches and generated tsunami. 
 \newblock {\em J. Comput. Phys.} {\bf 227}(16): 7720--7754,  2008.

\bibitem{GraCui}
J.~M. N.~T. Gray and X.~Cui.
\newblock Weak, strong and detached oblique shocks in gravity-driven granular
  free-surface flows.
\newblock {\em J. Fluid Mech.}, {\bf 579}:113--136, 2007.

\bibitem{GrTaNo}
J.~M. N.~T. Gray, Y.-C. Tai, and S.~Noelle.
\newblock Shock waves, dead zones and particle-free regions in rapid granular
  free-surface flows.
\newblock {\em J. Fluid Mech.}, {\bf 491}:161--181, 2003.

\bibitem{Gw1} P.~Gwiazda.
\newblock An existence result for a model of granular material with non-constant density. 
\newblock {\em Asymptot. Anal.}, {\bf 30}(1):  43--60, 2002.

\bibitem{Gw2}  P.~Gwiazda.
\newblock  On measure-valued solutions to a two-dimensional gravity-driven avalanche flow model. 
\newblock {\em Math. Methods Appl. Sci. } {\bf 28}(18): 2201--2223, 2005.
 
  

\bibitem{HuWaPu}
K.~Hutter,  Y.~Wang,  S.~P.~Pudasaini. 
\newblock The Savage-Hutter avalanche model: how far can it be pushed? 
\newblock {\em Philos. Trans. R. Soc. Lond. Ser. A Math. Phys. Eng. Sci.} {\bf 363}(1832):  1507--1528, 2005. 

\bibitem{JuMuGaNa} C.~Juez, J.~Murillo, P.~Garc'a-Navarro. 
\newblock 2D simulation of granular flow over irregular steep slopes using global and local coordinates. 
\newblock {\em J. Comput. Phys.} {\bf 255}:166--204, 2013.


\bibitem{PeBoMa}
M.~Pelanti, F.~Bouchut, A.~Mangeney.
\newblock A Roe-type scheme for two-phase shallow granular flows over variable topography. 
\newblock {\em M2AN Math. Model. Numer. Anal.} {\bf 42}(5): 851--885,  2008.

\bibitem{Sch}
V.~Scheffer.
\newblock An inviscid flow with compact support in space-time. 
\newblock {\em J. Geom. Anal.} 3(4):  343--401, 1993.

\bibitem{Shn}
A.~Shnirelman. 
\newblock On the nonuniqueness of weak solution of the Euler equation. 
\newblock {\em Comm. Pure Appl. Math. } {\bf 50}(12): 1261--1286, 1997.


\bibitem{Wi}
E.~Wiedemann
\newblock Existence of weak solutions for the incompressible Euler equations. 
\newblock {\em Ann. Inst. H. Poincar\'e Anal. Non Lin\'eaire} {\bf 28}(5): 727--730, 2011.

\bibitem{ZaPeTaNi}
N.~Zahibo, E.~Pelinovsky, T.~Talipova, and I.~Nikolkin.
\newblock {S}avage--{H}utter model for avalanche dynamics in inclined channels:
  {A}nalytical solutions.
\newblock {\em J. Geophys. Res.}, {\bf 115}:B3402, 1--18, 2010.

\end{thebibliography}
%\bibliographystyle{plain}

\end{document}